\documentclass[12pt]{article}
\usepackage[latin1]{inputenc}
\usepackage{amsmath, amsthm, amssymb}
\begin{document}

\bibliographystyle{plain}
\renewcommand{\Im}{\mathop{\rm Im }}
\renewcommand{\Re}{\mathop{\rm Re }}
\newcommand{\ra}{\mathop{\rightarrow }}
\newcommand{\supp}{\mathop{\rm supp}}
\newcommand{\sgn}{\mathop{\rm sgn }}
\newcommand{\card}{\mathop{\rm card }}
\newcommand{\KM}{\mbox{\rm KM}}
\newcommand{\diam}{\mathop{\rm diam}}
\newcommand{\diag}{\mathop{\rm diag}}
\newcommand{\tr}{\mathop{\rm tr}}
\newcommand{\Tr}{\mathop{\rm Tr}}
\newcommand{\dd}{ {\rm d} }
\newcommand{\id}{\mbox{\rm1\hspace{-.2ex}\rule{.1ex}{1.44ex}}
   \hspace{-.82ex}\rule[-.01ex]{1.07ex}{.1ex}\hspace{.2ex}}
\renewcommand{\P}{\mathop{\rm Prob}}
\newcommand{\V}{\mathop{\rm Var}}
\newcommand{\cps}{{\stackrel{{\rm p.s.}}{\longrightarrow}}}
\newcommand{\limm}{\mathop{\rm l.i.m.}}
\newcommand{\cloi}{{\stackrel{{\rm loi}}{\rightarrow}}}
\newcommand{\bra}{\langle\,}
\newcommand{\ket}{\,\rangle}
\newcommand{\ind}{{\bf 1}}
\newcommand{\obl}{/\!/}
\newcommand{\mapdown}[1]{\vbox{\vskip 4.25pt\hbox{\bigg\downarrow
  \rlap{$\vcenter{\hbox{$#1$}}$}}\vskip 1pt}}
\newcommand{\tab}{&\!\!\!}

\newcommand{\tabb}{&\!\!\!\!\!}
\renewcommand{\d}{\displaystyle}
\newcommand{\epreuve}{\hspace{\fill}$\bigtriangleup$}
\newcommand{\demo}{{\par\noindent{\em D\'emonstration~:~}}}
\newcommand{\solu}{{\par\noindent{\em Solution~:~}}}
\newcommand{\NB}{{\par\noindent{\bf Remarque~:~}}}
\newcommand{\const}{{\rm const}}
\newcommand{\cA}{{\cal A}}
\newcommand{\cB}{{\cal B}}
\newcommand{\cC}{{\cal C}}
\newcommand{\cD}{{\cal D}}
\newcommand{\cE}{{\cal E}}
\newcommand{\cF}{{\cal F}}
\newcommand{\cG}{{\cal G}}
\newcommand{\cH}{{\cal H}}
\newcommand{\cI}{{\cal I}}
\newcommand{\cJ}{{\cal J}}
\newcommand{\cK}{{\cal K}}
\newcommand{\cL}{{\cal L}}
\newcommand{\cM}{{\cal M}}
\newcommand{\cN}{{\cal N}}
\newcommand{\cO}{{\cal O}}
\newcommand{\cP}{{\cal P}}
\newcommand{\cQ}{{\cal Q}}
\newcommand{\cR}{{\cal R}}
\newcommand{\cS}{{\cal S}}
\newcommand{\cT}{{\cal T}}
\newcommand{\cU}{{\cal U}}
\newcommand{\cV}{{\cal V}}
\newcommand{\cW}{{\cal W}}
\newcommand{\cX}{{\cal X}}
\newcommand{\cY}{{\cal Y}}
\newcommand{\cZ}{{\cal Z}}
\newcommand{\bA}{{\bf A}}
\newcommand{\bB}{{\bf B}}
\newcommand{\bC}{{\bf C}}
\newcommand{\bD}{{\bf D}}
\newcommand{\bE}{{\bf E}}
\newcommand{\bF}{{\bf F}}
\newcommand{\bG}{{\bf G}}
\newcommand{\bH}{{\bf H}}
\newcommand{\bI}{{\bf I}}
\newcommand{\bJ}{{\bf J}}
\newcommand{\bK}{{\bf K}}
\newcommand{\bL}{{\bf L}}
\newcommand{\bM}{{\bf M}}
\newcommand{\bN}{{\bf N}}
\newcommand{\bP}{{\bf P}}
\newcommand{\bQ}{{\bf Q}}
\newcommand{\bR}{{\bf R}}
\newcommand{\bS}{{\bf S}}
\newcommand{\bT}{{\bf T}}
\newcommand{\bU}{{\bf U}}
\newcommand{\bV}{{\bf V}}
\newcommand{\bW}{{\bf W}}
\newcommand{\bX}{{\bf X}}
\newcommand{\bY}{{\bf Y}}
\newcommand{\bZ}{{\bf Z}}
\newcommand{\bu}{{\bf u}}
\newcommand{\bv}{{\bf v}}
\newfont{\msbm}{msbm10 scaled\magstep1}
\newfont{\msbms}{msbm7 scaled\magstep1} 
\newcommand{\bbA}{\mbox{$\mbox{\msbm A}$}}
\newcommand{\bbB}{\mbox{$\mbox{\msbm B}$}}
\newcommand{\bbC}{\mbox{$\mbox{\msbm C}$}}
\newcommand{\bbD}{\mbox{$\mbox{\msbm D}$}}
\newcommand{\bbE}{\mbox{$\mbox{\msbm E}$}}
\newcommand{\bbF}{\mbox{$\mbox{\msbm F}$}}
\newcommand{\bbG}{\mbox{$\mbox{\msbm G}$}}
\newcommand{\bbH}{\mbox{$\mbox{\msbm H}$}}
\newcommand{\bbI}{\mbox{$\mbox{\msbm I}$}}
\newcommand{\bbJ}{\mbox{$\mbox{\msbm J}$}}
\newcommand{\bbK}{\mbox{$\mbox{\msbm K}$}}
\newcommand{\bbL}{\mbox{$\mbox{\msbm L}$}}
\newcommand{\bbM}{\mbox{$\mbox{\msbm M}$}}
\newcommand{\bbN}{\mbox{$\mbox{\msbm N}$}}
\newcommand{\bbO}{\mbox{$\mbox{\msbm O}$}}
\newcommand{\bbP}{\mbox{$\mbox{\msbm P}$}}
\newcommand{\bbQ}{\mbox{$\mbox{\msbm Q}$}}
\newcommand{\bbR}{\mbox{$\mbox{\msbm R}$}}
\newcommand{\bbS}{\mbox{$\mbox{\msbm S}$}}
\newcommand{\bbT}{\mbox{$\mbox{\msbm T}$}}
\newcommand{\bbU}{\mbox{$\mbox{\msbm U}$}}
\newcommand{\bbV}{\mbox{$\mbox{\msbm V}$}}
\newcommand{\bbW}{\mbox{$\mbox{\msbm W}$}}
\newcommand{\bbX}{\mbox{$\mbox{\msbm X}$}}
\newcommand{\bbY}{\mbox{$\mbox{\msbm Y}$}}
\newcommand{\bbZ}{\mbox{$\mbox{\msbm Z}$}}
\newcommand{\bbsA}{\mbox{$\mbox{\msbms A}$}}
\newcommand{\bbsB}{\mbox{$\mbox{\msbms B}$}}
\newcommand{\bbsC}{\mbox{$\mbox{\msbms C}$}}
\newcommand{\bbsD}{\mbox{$\mbox{\msbms D}$}}
\newcommand{\bbsE}{\mbox{$\mbox{\msbms E}$}}
\newcommand{\bbsF}{\mbox{$\mbox{\msbms F}$}}
\newcommand{\bbsG}{\mbox{$\mbox{\msbms G}$}}
\newcommand{\bbsH}{\mbox{$\mbox{\msbms H}$}}
\newcommand{\bbsI}{\mbox{$\mbox{\msbms I}$}}
\newcommand{\bbsJ}{\mbox{$\mbox{\msbms J}$}}
\newcommand{\bbsK}{\mbox{$\mbox{\msbms K}$}}
\newcommand{\bbsL}{\mbox{$\mbox{\msbms L}$}}
\newcommand{\bbsM}{\mbox{$\mbox{\msbms M}$}}
\newcommand{\bbsN}{\mbox{$\mbox{\msbms N}$}}
\newcommand{\bbsO}{\mbox{$\mbox{\msbms O}$}}
\newcommand{\bbsP}{\mbox{$\mbox{\msbms P}$}}
\newcommand{\bbsQ}{\mbox{$\mbox{\msbms Q}$}}
\newcommand{\bbsR}{\mbox{$\mbox{\msbms R}$}}
\newcommand{\bbsS}{\mbox{$\mbox{\msbms S}$}}
\newcommand{\bbsT}{\mbox{$\mbox{\msbms T}$}}
\newcommand{\bbsU}{\mbox{$\mbox{\msbms U}$}}
\newcommand{\bbsV}{\mbox{$\mbox{\msbms V}$}}
\newcommand{\bbsW}{\mbox{$\mbox{\msbms W}$}}
\newcommand{\bbsX}{\mbox{$\mbox{\msbms X}$}}
\newcommand{\bbsY}{\mbox{$\mbox{\msbms Y}$}}
\newcommand{\bbsZ}{\mbox{$\mbox{\msbms Z}$}}

%
\def\eurtoday{\number\day \space\ifcase\month\or
 January\or February\or March\or April\or May\or June\or
 July\or August\or September\or October\or November\or December\fi
 \space\number\year}
%
%
\def\aujourdhui{\number\day \space\ifcase\month\or
 janvier\or f{\'e}vrier\or mars\or avril\or mai\or juin\or
 juillet\or ao{\^u}t\or septembre\or octobre\or novembre\or d{\'e}cembre\fi
 \space\number\year}
\textheight=21cm
\textwidth=16cm 
\voffset=-1cm
\newtheorem{theo}{Theorem}[section]
\newtheorem{pr}{Proposition}[section]
\newtheorem{cor}{Corollary}[section]
\newtheorem{lem}{Lemma}[section]
\newtheorem{defn}{Definition}[section]
\hoffset=-1,5cm
\parskip=4mm
\title{A functional approach for random walks in random sceneries}
\author{C. Dombry  and N. Guillotin-Plantard \protect\hspace{1cm}}
\date{}
\maketitle
~\\
~\\
~\\
~\\ 
Authors' addresses:\\
C. Dombry (corresponding author), Université de Poitiers, Laboratoire de Mathématiques et Applications, Téléport 2 - BP 30179,
Boulevard Marie et Pierre Curie, 86962 Futuroscope Chasseneuil Cedex, France, \\*
e-mail: clement.dombry@math.univ-poitiers.fr\\
N. Guillotin-Plantard, Universit\'e de Lyon, Institut Camille Jordan, 43 bld du 11 novembre 1918, 69622 Villeurbanne, France, e-mail: nadine.guillotin@univ-lyon1.fr \\

\begin{abstract}
A functional approach for the study of the random walks in random sceneries (RWRS) is proposed. Under fairly general assumptions on the random walk and on the random scenery, functional limit theorems are proved. The method allows to study separately the convergence of the walk and of the scenery: on the one hand, a general criterion for the convergence of the local time of the walk is provided, on the other hand, the convergence of the random measures associated with the scenery is studied. This functional approach is robust enough to recover many of the known results on RWRS as well as new ones, including the case of many walkers evolving in the same scenery. 
\end{abstract}

\section{Introduction}
Let $S=(S_n)_{n\in\bbsN}$ be a {\it random walk} in $\bbZ^d$ starting from zero. The term "random walk" has to be understood in the broad sense of a discrete process, not necessarily with independent and identically distributed (i.i.d.) increments. Let $\xi=\{\xi_x, x\in\bbZ^d\}$ be a collection of random variables called the {\it random scenery}. The {\it random walk in random scenery} (RWRS) is the cumulative process
$$Z_n=\sum_{k=0}^n \xi_{S_k}.$$ 
Motivated by the construction of a new class of self-similar stationary increments processes, Kesten and Spitzer \cite{KS} and Borodin \cite{Bor1,Bor2} introduced RWRS in dimension one and proved functional limit theorems. More precisely, they considered the case when the random walk has i.i.d. increments belonging to the domain of attraction of a stable distribution with index $\alpha\in (1,2]$, the random scenery is i.i.d. and belongs to the domain of attraction of a stable distribution with index $\beta\in (0,2]$; the scenery and the walk are independent. Under these assumptions, the renormalised process $(n^{-\delta}Z_{[nt]})_{t\geq 0}$ is shown to converge to a $\delta$-self-similar stationary increments process where $\delta=1-\alpha^{-1}+(\alpha\beta)^{-1}$.

Since then, many authors considered  RWRS models and related functional limit theorems. Bolthausen \cite{bol} studied the case of a recurrent $\bbZ^2$-random walk. Maejima \cite{Mae} considered the case of a random walk in dimension one and of a multivariate random scenery in link with operator self-similar random processes. All  these results require independence properties on the increments of the walk and on the scenery. This constraint can be relaxed by assuming some reasonable dependence structure: in Wang \cite{WA}, the increments of the random walk are assumed to be strongly correlated, Gaussian and in the domain of attraction of a fractional Brownian motion; it evolves in a i.i.d. square integrable random scenery. In Lang \& Xahn \cite{LX} a functional limit theorem is proved for an i.i.d. increments $\bbZ$-random walk in a strongly correlated scenery satisfying a non central limit theorem of Dobrushin \& Major \cite{DM}. Weakly dependent sceneries have also recently been studied by Guillotin-Plantard \& Prieur \cite{GPP1,GPP2} under a $\theta$-mixing condition. Let us also mention the work of Le Borgne \cite{Leb} in which the random scenery verifies some strong decorrelation condition.

One dimensional RWRS arise in the study of random walks evolving on oriented versions of $\bbZ^2$ (Guillotin-Plantard \& Le Ny \cite{GPLN1,GPLN2}; Pène \cite{PE2}) and in the context of charged polymers (Chen \& Khoshnevisan \cite{CHK}). Motivated by models for traffic in a network, Cohen \& Samorodnitsky \cite{CS} introduced a random reward schema consisting of sums of independent copies of RWRS and studied its convergence. Extensions of this work have been considered in Guillotin-Plantard \& Dombry \cite{domgui} and Cohen \& Dombry \cite{CD}.

In this paper, a functional approach for RWRS is developed and the convergence of the finite-dimensional distributions of RWRS under very general assumptions that cover many of the above results are proved. A functional approach using stochastic integration was suggested  by Borodin \cite{Bor2} and Cadre \cite{C} in the case when the random scenery is a martingale. Our approach is rather based on the convergence of random measures to random noises (see \cite{DOM}) and is more adapted when the scenery is not a martingale (including weakly and strongly correlated  RS). Under the assumption of independence of the RW and of the RS (assumed in all the above cited works), it allows to study separately the convergence of the local time associated with the random walk on the one hand, and the convergence of the random measures associated with the random scenery on the other hand. The technique is also robust enough to deal with a variant of RWRS: the case when many walkers evolve in a single scenery, which might appear as a more realistic model of the random reward scheme considered in \cite{CS} since all the users share the same network. To our knowledge, these are the first results in that direction.
  
Our paper is organized as follows. In Section 2, after introducing a general criterion for the $L^p$-convergence of the local time of the random walk, we verify that it applies in the classical cases above-mentioned. In Section 3 a sequence of random measures is associated with the random scenery and conditions to ensure its convergence to a stable, Gaussian or fractional Gaussian random noise are given. In Section 4 we explain how these results allow to recover limit theorems for RWRS and also consider the case of many walkers evolving in one single random scenery.

\section{Convergence in $L^p$ of the local times of the random walk}
\subsection{A general criterion}
Let $S=(S_n)_{n\in \bbsN}$ be a $\bbZ$-valued random sequence. We define the discrete local time at time $n\in\bbN$ and point $x\in\bbZ$ by
$$N(n,x)={\rm card}\{k\in [\!| 0,n|\!]; S_k=x\}.$$
\begin{pr}\label{gencrit}
Suppose that the following assumption $\mbox{\rm\bf(RW)}$ is satisfied:
\begin{itemize}
\item{\bf(RW1)} There is some sequence $a_n$ verifiying $a_n\to\infty$ and $n^{-1}a_n\to 0$ such that the renormalized process $\left(a_n^{-1}S_{[nt]}\right)_{t\geq 0}$ converges in $\cD([0,\infty))$ to
some process $(Y_t)_{t\geq 0}$ admitting a local time $(L(t,x))_{t\geq 0, x\in\bbsR}$.
\item{\bf(RW2)} There is some $p\geq 1$ such that for all $M>0$:
$$\lim_{\delta\to 0} \limsup_{n\to\infty} \int_{[-M,M]} \bbE |L_n(t,x)-L_n(t,[x]_\delta)|^p dx=0$$
where $L_n(t,x)$ is the rescaled discrete local time
\begin{equation}\label{eq:LTscaling}
L_n(t,x)= n^{-1}a_n N([nt],[a_n x])
\end{equation}
and $[x]_\delta=\delta [\delta^{-1}x]$.
\end{itemize}
Then, for every $m\geq 1$, for every $\theta_1\in\bbR,\ldots, \theta_m\in\bbR$, for every $(t_1,\ldots,t_m)$ such that $0<t_1<\ldots<t_m$, the following convergence holds in $L^p, p\geq 1$: 
$$\sum_{i=1}^m \theta_i L_n(t_i,.)\stackrel{\cL}{\Longrightarrow} \sum_{i=1}^m \theta_i L(t_i,.) \quad {\mbox as}\quad n\to\infty.$$ 
\end{pr}

\noindent{\bf Remark:} A direct application of Lebesgue's dominated convergence Theorem shows that Assumption ${\bf (RW2)}$ is a consequence of the following two conditions: \\
${\bf(RW2.a)}$ There is some $A>0$ such that for all $n\geq 1$, 
$$\sup_{x\in\bbsR} \bbE |L_n(t,x)|^p \leq A,$$
${\bf(RW2.b)}$ For every $t>0$, for almost every $x\in\bbR$, 
$$\lim_{y\to 0} \limsup_{n\to\infty}  \bbE |L_n(t,x)-L_n(t,x+y)|^p=0.$$
We notice that if these conditions are verified for some $p>1$ then they hold for any $p'\in [1, p]$.\\*

\noindent{\bf Proof:}\\
We prove the theorem for $m=1$ and $t_1=t>0$. Generalizations to any $m\geq 1$ can be easily obtained.\\* 
We denote by $q$ the conjugate of $p$.\\*
According to Theorem 2.4 of De Acosta (1970) we need to check two things:\\
- first that there is a $w^\star$-dense set $D\in L^q(\bbR)$ such that for each $f\in D$,
$$ \int_{\bbsR} f(x)L_n(t,x)dx \rightarrow \int_{\bbsR} f(x)L(t,x)dx,$$
- second that the sequence $\left(L_n(t,.)\right)_{n\geq 1}$ is flatly concentrated, that is, for every $\varepsilon>0$, there is a finite dimensional subspace $F$ of $L^p(\bbR)$ such that
$$\liminf_{n\to\infty} \bbP\left[L_n(t,.)\in F^\varepsilon \right]\geq 1-\varepsilon $$
where $F^\varepsilon$ is the $\varepsilon$-neighborhood of $F$.

The first assertion follows from the invariance principle ${\bf(RW1)}$ for $a_n^{-1}S_{[n.]}$. We take $D$ as the set of continuous functions with bounded support and observe that the definition of the discrete occupation time $L_n$ and the uniform continuity of $f$ imply
\begin{eqnarray*}
\int_{\bbsR} f(x)L_n(t,x)dx&=&n^{-1}\int_{\bbsR} f(a_n^{-1}x)N([nt],[x])dx \\
&=& n^{-1}\int_{\bbsR} f(a_n^{-1}[x])N([nt],[x])dx +o(1) \\
&=& n^{-1}\sum_{k=1}^{[nt]} f(a_n^{-1}S_k) +o(1)\\
&=& \int_0^t  f(a_n^{-1}S_{[nu]})du +o(1)
\end{eqnarray*}
From the continuity of the application 
$$\cD([0,t])\rightarrow \bbR, \phi\mapsto \int_0^ t f\circ \phi,$$
from the invariance principle ${\bf (RW1)}$ and the existence of the local time $L$ for $Y$,  we deduce
$$\int_{\bbsR} f(x)L_n(t,x)dx \rightarrow \int_0^t  f(Y_u)du=\int_{\bbsR}f(x)L(t,x)dx.$$ 

We turn to the second assertion and consider the partition $x_j=j\delta$, $j\in\bbZ$ for a given $\delta>0$. 
Define the  finite dimensional space $F$ of $L^p(\bbR)$ as the space of functions that are constant on any interval $(x_j,x_{j+1}), j\in\bbZ$ and that vanish outside $[-M,M]$ (for some $M>0$ to be specified later).
Define $\tilde L_n$ by $\tilde L_n(x)=L_n(t,x_j)$ if $x_j\leq x < x_{j+1}$, i.e. $\tilde L_n(x)=L_n(t,[x]_\delta)$. 
Notice that $\tilde L_n\in F$ on the event 
$$\left\{L_n(t,.)\equiv 0 \mbox{\ on\ } \bbR\setminus[-M,M]\right\}= \left\{\max_{0\leq k\leq [nt]} |a_n^{-1}S_k| \leq M\right\}.$$
Therefore, 
\begin{eqnarray*}
& &\bbP\left[L_n(t,.)\notin F^\varepsilon \right]\\
&\leq& \bbP \left[ \max_{0\leq k\leq [nt]} |a_n^{-1}S_k| > M \right] + \bbP \left[ \max_{0\leq k\leq [nt]} |a_n^{-1}S_k| \leq M\ \mbox{and}\ |\!|L_n(t,.)-\tilde L_n|\!|_{L^p} \geq \varepsilon \right]  
\end{eqnarray*}
By the invariance principle ${\bf (RW1)}$, 
$$\lim_{n\to\infty}\bbP \left[ \max_{0\leq k\leq [nt]} |a_n^{-1}S_k| > M \right]=\bbP \left[ \max_{0\leq s\leq t} |Y_s| > M \right]$$
is small for large $M$. Hence it remains only to show that
$$\lim_{\delta\to 0} \limsup_{n\to\infty} \bbP \left[ \max_{0\leq k\leq [nt]} |a_n^{-1}S_k| \leq M\ \mbox{and}\ |\!|L_n(t,.)-\tilde L_n|\!|_{L^p} \geq \varepsilon \right]=0.$$
Using Markov's inequality, we get
\begin{eqnarray*}
& &\bbP \left[ \max_{0\leq k\leq [nt]} |a_n^{-1}S_k| \leq M\ \mbox{and}\ |\!|L_n(t,.)-\tilde L_n|\!|_{L^p} \geq \varepsilon \right]\\
&\leq& \varepsilon^{-p} \bbE\left[ \int_{[-M,M]} |L_n(t,x)-\tilde L_n(t,[x]_\delta)|^p dx \right]\\
\end{eqnarray*}
and the result follows from Assumption ${\bf(RW2)}$. \ \ $\Box$

\subsection{Local time of a random walk in the domain of attraction of a stable Lévy motion}
Let $S=(S_n)_{n\in\bbsN}$ be a {\bbZ}-random walk such that 
$$\left\{\begin{array}{l}S_0=0, \\ 
S_n=\displaystyle\sum_{k=1}^n X_k\ ,\ n\geq 1,
\end{array}\right.$$ 
where the $X_i$ are i.i.d. integer-valued random variables belonging to the normal domain of attraction of a strictly stable distribution $\cS_{\alpha}(\sigma,\nu,0)$ with parameters 
$$\alpha\in (1,2], \sigma>0 \ \mbox{and}\ \nu\in[-1,1].$$
This means that the following weak convergence holds:
\begin{equation}\label{eq1.04}
n^{-\frac{1}{\alpha}} S_n \mathop{\Longrightarrow}_{n\rightarrow\infty}^{\cL} Z_\alpha, 
\end{equation}
where the characteristic function of the random variable $Z_{\alpha}$ is given by
\begin{equation}\label{eq1.05}
\forall u\in\bbR,\ \  \bbE\exp(iuZ_\alpha)=\exp\left(-\sigma^\alpha|u|^\alpha(1-i\nu\tan(\frac{\pi\alpha}{2})\sgn(u))\right). 
\end{equation}
Assumption ${\bf (RW1)}$ holds with $a_n=n^{1/\alpha}$ and $(Y(t))_{t\geq 0}$ an $\alpha$-stable Lévy process such that $Y(0)=0$ and $Y(1)$ is distributed according to $Z_\alpha$.\\*
The local time $(N(n,x))_{n\in\bbsN;x\in\bbsZ}$ of the random walk $S$ is defined as in the previous section.
Its rescaled local time $(L_n(t,x))_{t\geq 0,x\in\bbsR}$ is defined as
$$L_n(t,x)= n^{-(1-1/\alpha)} N([nt],[n^{1/\alpha} x]).$$
We denote by $(L(t,x))_{t\geq 0,x\in\bbsR}$ the jointly continuous version of the local time of the process $Y$. We are interested in the convergence in $L^p, p\geq 1$ of the rescaled local time 
$L_n$ to $L$. By applying the general criterion of Section 2.1, we get the following result.
\begin{pr}\label{stalev}
For every $m\geq 1$, for every $\theta_1\in\bbR,\ldots, \theta_m\in\bbR$, for every $(t_1,\ldots,t_m)$ such that $0<t_1<\ldots<t_m$, the following convergence holds in $L^p, p\in [1,2]$: 
$$\sum_{i=1}^m \theta_i L_n(t_i,.)\stackrel{\cL}{\Longrightarrow} \sum_{i=1}^m \theta_i L(t_i,.) \quad {\mbox as}\quad n\to\infty.$$ 
\end{pr}
\noindent{\bf Proof:}\\*
From the remark following Proposition \ref{gencrit}, it is enough to prove that conditions ${\bf(RW2.a-b)}$ are satisfied in our setting with $p=2$. \\*
\noindent{\bf Condition (RW2.a):} Kesten and Spitzer proved in \cite{KS} ((2.10) in Lemma 1) that if $\tau_{x}$ denotes the hitting time of the point $x$ by the random walk $S$, then the following inequality holds for any $x\in\bbZ$, $r\geq 0$ and $n\geq 1$,
$$\bbP(N(n,x)\geq r)\leq \bbP(N(n+1,0)\geq r) \bbP(\tau_x\leq n+1).$$
It implies that the moment
 of order 2 of the random variable $N(n,x)$ is uniformly bounded by the moment of order 2 of $N(n+1,0)$ which is equivalent to $Cn^{2(1-1/\alpha)}$ for some $C>0$ (see (2.12) in Lemma 1 of \cite{KS}). 
\\*
\noindent{\bf Condition (RW2.b):} We know (see Lemma 2 and 3 in \cite{KS}) that there exists some constant $C>0$ independent of $x,y\in \bbZ$ and $n\geq 1$ s.t.
$$\bbE(|N(n,[n^{1/\alpha}x])-N(n,[n^{1/\alpha}(x+y)])|^2)\leq C \left|[n^{1/\alpha}x]-[n^{1/\alpha}(x+y)]\right|^{\alpha-1} n^{1-1/\alpha}.$$
Then, we have
\begin{eqnarray*}
& &\lim_{y\to 0} \limsup_{n\to\infty} \bbE(|L_n(t,x)-L_n(t,x+y)|^2)\\
&\leq &C \lim_{y\to 0}\limsup_{n\to\infty} 
\left|\frac{[n^{1/\alpha}x]}{n^{1/\alpha}}-\frac{[n^{1/\alpha}(x+y)]}{n^{1/\alpha}}\right|^{\alpha-1}\\
&=& C \lim_{y\to 0} |y|^{\alpha-1}=0\ \ (\alpha>1\ \ \mbox{\rm by hypothesis}).\ \ \  \Box
\end{eqnarray*}

\subsection{Local time of a strongly correlated walk in the domain of attraction of a fractional Brownian motion}
Let $S=(S_n)_{n\in\bbsN}$ be a {\bbZ}-random walk such that 
$$\left\{\begin{array}{l}S_0=0, \\ 
S_n=\left[\displaystyle\sum_{k=1}^n X_k\right]\ ,\ n\geq 1,
\end{array}\right.$$ 
where $(X_k)_{k\geq 1}$ is a stationary Gaussian sequence with mean 0 and correlations $r(i-j)=\bbE(X_i X_j)$
satisfying as $n\to\infty$,
$$\sum_{i=1}^n\sum_{j=1}^n r(i-j)\sim n^{2H} l(n)$$
with $0<H<1$ and $l$ is a slowly varying function at $\infty$. Assumption ${\bf(RW1)}$ is satisfied with $a_n=n^{H} l(n)^{1/2}$ and $Y$ the standard fractional Brownian motion $(B_{H}(t))_{t\geq 0}$ of index $0<H<1$, i.e. the centered Gaussian process such that $B_H(0)=0$ and $\bbE(|B_{H}(t_2)-B_{H}(t_1)|^2)=|t_2-t_1|^{2H}$. 
The local time $N_H(n,x)$ with $n\in\bbN$ and $x\in\bbZ$ is defined as the number of visits of $(S_n)_{n\in\bbsN}$ at point $x$ up to time n.
Its rescaled local time $(L_{n,H}(t,x))_{t\geq 0,x\in\bbsR}$ is defined as
$$L_{n,H}(t,x)= n^{-(1-H)} l(n)^{1/2} N([nt],[n^{H}l(n)^{1/2} x]).$$
We denote by $(L_H(t,x))_{t\geq 0,x\in\bbsR}$ the jointly continuous version of the local time of the process $B_{H}$. 
We are interested in the convergence in $L^p, p\geq 1$ of the rescaled local time 
$L_{n,H}$ to $L_H$. By applying the general criterion of Section 2.1, we get the following result.
\begin{pr}\label{Gauss}
For every $m\geq 1$, for every $\theta_1\in\bbR,\ldots, \theta_m\in\bbR$, for every $(t_1,\ldots,t_m)$ such that $0<t_1<\ldots<t_m$, the following convergence holds in $L^p, p\in [1,2]$: 
$$\sum_{i=1}^m \theta_i L_{n,H}(t_i,.) \stackrel{\cL}{\Longrightarrow}\sum_{i=1}^m \theta_i L_H(t_i,.) \quad {\mbox as}\quad n\to\infty.$$ 
\end{pr}
\noindent{\bf Proof:}\\*
Conditions $\bf{(RW2.a-b)}$ directly follow from Lemma 4.4 and Lemma 4.6 of Wang \cite{WA}. $\Box$

\subsection{The local time of many independent random walks}
We consider independent and identically distributed $\bbZ-$random walks 
$$S^{(i)}=(S_n^{(i)})_{n\in\bbsN},\quad i\geq 1,$$
verifying assumption {\bf({RW})}. The local time of the $i$-th random walk is denoted by $(N^{(i)}(n,x))_{n\in\bbsN;x\in\bbsZ}$ and the rescaled one by $(L_n^{(i)}(t,x))_{t\geq 0;x\in\bbsR}$.  The independent limit processes of the correctly renormalized random walks $S^{(i)}$ have local time denoted by $(L^{(i)}(t,x))_{t\geq 0;x\in\bbsR}$. We prove that
\begin{pr}\label{manywa}
Let $c_n\to\infty$ as $n\to \infty$. Suppose that the random walks $S^{(i)},i\geq 1$ are i.i.d. and satisfy assumption {\bf({RW})} for some $p\geq 1$. Suppose furthermore that for any fixed $t\geq 0$, 
\begin{equation}\label{eq:integLT}
\bbE\left[ \left(\int_{\bbsR} L^{(1)}(t,x)^pdx \right)^{1/p}\right]<\infty.
\end{equation}
Then, for every $m\geq 1$, for every $\theta_1\in\bbR,\ldots, \theta_m\in\bbR$, for every $(t_1,\ldots,t_m)$ such that $0<t_1<\ldots<t_m$, the following convergence holds in $L^p$: 
\begin{equation}\label{eq1}
\frac{1}{c_n}\sum_{j=1}^m \theta_j \sum_{i=1}^{c_n} L_n^{(i)}(t_j,.)\stackrel{{\cal L}}{\Longrightarrow} \sum_{j=1}^m \theta_j \bbE(L^{(1)}(t_j,.)) \quad {\mbox as}\quad n\to\infty.
\end{equation}
\end{pr}
\noindent{\bf Proof:}\\*
From Proposition \ref{gencrit} and Skorohod's representation theorem, there exists a probability space on which are defined  independent copies of the processes $L_n^{(i)}$ and $L^{(i)}$ denoted by $\tilde{L}_n^{(i)}$ and $\tilde{L}^{(i)}$ such that for every $i\geq 1$, the sequence 
$\sum_{j=1}^m\theta_j\tilde{L}_n^{(i)}(t_j,.)$ converges almost surely in $L^{p}$ to $\sum_{j=1}^m\theta_j\tilde{L}^{(i)}(t_j,.)$ as $n$ tends to infinity. So it is enough to prove the result for these later sequences.  \\*
{\bf Step 1:}\\* 
First, from triangular inequality, we have
\begin{eqnarray*}
& &\left|\!\left|\frac{1}{c_n}\sum_{j=1}^m \theta_j \sum_{i=1}^{c_n} \left[\tilde{L}_n^{(i)}(t_j,.)-\tilde{L}^{(i)}(t_j,.)\right] \right|\!\right|_{L^p(\bbsR)}\\
&\leq & \frac{1}{c_n}\sum_{i=1}^{c_n} |\!|\sum_{j=1}^m \theta_j\left(\tilde{L}_n^{(i)}(t_j,.)-\tilde{L}^{(i)}(t_j,.)\right)|\!|_{L^p(\bbsR)}.
\end{eqnarray*}
and we prove that this quantity converges to zero as $n\to\infty$. To see this, let 
$$X_n^{(i)}=|\!|\sum_{j=1}^m \theta_j\left(\tilde{L}_n^{(i)}(t_j,.)-\tilde{L}^{(i)}(t_j,.)\right)|\!|_{L^p(\bbsR)}.$$
The triangular array of random variables $\{X_n^{(i)}, n\geq 1, i\in \{1,\ldots,c_n\}\}$ is such that:\\
{}- on each row (fixed $n$), the random variables $X_n^{(i)}, i\in\{1,\ldots,c_n\}$ are i.i.d..\\
{}- on each column (fixed $i$), $(X_n^{(i)})_{n\geq 1}$ weakly converges to zero as $n\to\infty$.\\
Under these general assumptions, the following weak law of large numbers holds 
\begin{equation}\label{eq:WLLN}
\frac{1}{c_n}\sum_{i=1}^{c_n} X_n^{(i)}\stackrel{{\cal L}}{\Longrightarrow} 0, 
\end{equation}
implying the required result for the local time.\\
We finally give some indications for the proof of the weak law of large numbers. Denote by $A_n$ the random variable in the left hand side of (\ref{eq:WLLN}). Using the independence of the $X_n^{(i)}$'s, its characteristic function is given by
$$\bbE\left[\exp(iuA_n) \right]=\left(\bbE\left[\exp(iuc_n^{-1}X^{(1)}_n) \right] \right)^{c_n}$$
and hence
$$\left|\bbE\left[\exp(iuA_n) \right]-1\right|\leq \bbE\left[ \min\left( c_n|\exp(iuc_n^{-1}X^{(1)}_n)-1|,2 \right) \right] \leq \bbE\left[ \min\left(|u|X^{(1)}_n,2 \right) \right] .$$
The function $x\rightarrow\min\left(|u| x,2 \right)$ is continuous and bounded on $\bbR^{+}$. Then, since the sequence  $(X^{(1)}_n)_{n\geq 1}$ weakly converges to $0$, we deduce that $(\bbE\left[\exp(iuA_n) \right])_{n\geq 1}$ converges to $1$ as $n\to\infty$ and this proves the weak law of large numbers.\\
\noindent{\bf Step 2:}\\*
From Assumption (\ref{eq:integLT}) and the strong law of large numbers in the Banach space $L^p(\bbR)$ (see \cite{led} for instance), the sequence 
$$\frac{1}{c_n}\sum_{j=1}^m \theta_j \sum_{i=1}^{c_n}\tilde{L}^{(i)}(t_j,.)$$
converges almost surely in $L^p$ to the function
$$\sum_{j=1}^m \theta_j \bbE(\tilde{L}^{(1)}(t_j,.)).$$
\noindent {\bf Conclusion:}\\*
Combining steps 1 and 2 gives (\ref{eq1}). \ \ $\Box$

\noindent
{\bf Examples:}
The above proposition directly applies in the case of random walks with i.i.d. increments in the domain of attraction of a stable Lévy motion (see Section 2.2) or of a strongly correlated Gaussian random walk in the domain of attraction of a fractional Brownian motion (see Section 2.3). In both cases, equation (\ref{eq:integLT}) is satisfied for every $p\in [1,2]$:  Lemma 2.1 of \cite{domgui} states that the local time of a stable Lévy motion is $L^p-$integrable and  Theorem 3.1 of \cite{CS} gives a similar result for the fractional Brownian motion.

\section{Convergence of the random measures associated with the random scenery}
Let $\xi=\{\xi_x,x\in\bbZ\}$ be a family of real random variables. 
In the sequel, we consider $\mu_h=\mu_h(\xi)$ the random signed measure on $\bbR$ absolutely continuous with respect to Lebesgue measure with random density
\begin{equation}\label{eq:defmu}
\frac{d\mu_h}{dx}(x)=\gamma_hh^{-1}\sum_{k\in\bbsZ}\xi_k  \ind_{[hk,h(k+1))}(x), 
\end{equation}
where  $\gamma_h>0$ is a normalisation constant. \\
For a locally integrable function $f\in L^1_{loc}$, we want to define
\begin{equation}\label{eq:serie}
\mu_h[f]=\gamma_h\sum_{k\in\bbsZ}\xi_kh^{-1}\int_{hk}^{h(k+1)}f(x)dx.
\end{equation}
Denote by $\cF_{\mu_h}$ the set of functions $f\in L^1_{loc}$ such that this series is convergent in a sense that will be precised later according to the considered cases: either almost-sure semi-convergence or convergence in $L^2(\Omega)$.
Clearly, any integrable function with bounded support belongs to the set $\cF_{\mu_h}$.\\ 
The following scaling relation is worth noting:
\begin{equation}\label{eq:scaling}
\mu_h[f(c.)]=\frac{\gamma_h}{c\gamma_{ch}}\mu_{ch}[f(.)],
\end{equation}
whenever these quantities are well-defined.

We introduce the cumulative scenery $(w_x)_{x\in\bbsZ}$ 
$$w_x=\left\{\begin{array}{lll} 
\sum_{i=0}^{x-1}\xi_i  & \mbox{if\ } x >0 \\
\sum_{i=x}^{-1} \xi_k & \mbox{if\ } x< 0 \\
0 & \mbox{if\ } x= 0\\
 \end{array}\right.,\quad x\in\bbZ,$$
and we extend it to a continuous process by the linear interpolation
$$ w_x=w_{[x]}+(x-[x])(w_{[x]+1}-w_{[x]}).$$
Note that
\begin{equation}\label{eq:lien}
\mu_h[1_{[x_1, x_2]}]=\gamma_h (w_{h^{-1}x_2}-w_{h^{-1}x_1})=W_h(x_2)-W_h(x_1),\quad x_1,x_2\in\bbR, x_1\leq x_2.
\end{equation}
with $W_h$ the rescaled cumulative scenery $W_h(x)=\gamma_h w_{h^{-1}x}$.

\subsection{Independent and identically distributed scenery}
Let $1<\beta\leq 2$. We consider the case when the $\xi_x,x\in\bbZ$ satisfy the following assumption:
\ \\
\noindent
${\bf(RS1)}$  the $\xi_x$'s are i.i.d. random variables in the normal domain of attraction of the stable law $\cS_\beta(\sigma,\nu,0)$ with characteristic function given in (\ref{eq1.05}).\\*
\noindent In this case, with the scaling $\gamma_h=\sigma^{-1}h^{\frac{1}{\beta}}$, the rescaled cumulative scenery $W_h$ converges to a bilateral $\beta$-stable Lévy process $W$ such that $W(1)$ is distributed according to the stable distribution $\cS_\beta(1,\nu,0)$. In \cite{DOM}, the convergence of the random measures $\mu_h$ to a stable white noise was proved.
Here a function $f\in L^1_{loc}$ belongs to $\cF_{\mu_h}$ if the random series (\ref{eq:serie}) defining $\mu_h[f]$ converges almost surely.
\begin{pr}\label{clem}
Suppose that $\xi$ satisfies assumption ${\bf(RS1)}$ for some $1<\beta\leq 2$.\\
Then for any $h>0$, $L^\beta(\bbR)\subset \cF_{\mu_h}$  and if $f_n\to f$ in $L^\beta$ and $h_n\to 0$ as $n\to\infty$,
$$\mu_{h_n}[f_n]\stackrel{\cL}{\Longrightarrow} \int_{\bbsR}f(x)W(dx) $$
where $W$ denotes the independently scattered $\beta$-stable random noise on $\bbR$ with Lebesgue intensity and constant skewness $\nu$. 
\end{pr}

\subsection{A general criterion when the scenery is square integrable}
We consider the case of  a centered scenery with finite variance in the domain of attraction of some square integrable random process $W$. Our assumption $\mbox{\rm\bf (RS2)}$ is:
\begin{itemize}
\item{\bf(RS2.a)} The scenery $\xi$ is centered and square integrable.\\
\item{\bf(RS2.b)} The finite dimensional distributions of the rescaled cumulative scenery $W_h$ converge as $h\to 0$:\\
$$(W_h(x))_{x\in\bbsR} \stackrel{\cL}{\Longrightarrow} (W(x))_{x\in\bbsR}$$
where $W$ is a non-degenerate centered square integrable process.\\
\item{\bf(RS2.c)} There is some $C_1>0$ such that for any $h>0$, $K\geq 1$, $x_1\leq\cdots\leq x_{K+1}$, $\theta_1,\cdots,\theta_K\in\bbR$,  
$$ |\!|\sum_{i=1}^K \theta_i(W_h(x_{i+1})-W_h(x_i))|\!|_{L^2(\Omega)} \leq C_1 |\!| \sum_{i=1}^K \theta_i(W(x_{i+1})-W(x_i)) |\!|_{L^2(\Omega)}.$$
\end{itemize}
Note that $W(0)=0$ and that $W$ non degenerate means that for any pairwise distinct nonzero $x_1,\cdots,x_p\in\bbR$, the random variables $W(x_1),\cdots,W(x_p)$ are linearly independent in $L^2(\Omega)$.\\

We now define integration of deterministic functions with respect to $W$. This construction is classical in the case when $W$ is a Gaussian process (e.g. Brownian or fractional Brownian motion), and the reader can refer to \cite{Nua}, but we simply assume here that $W$ is centered square integrable and non degenerate.\\
Let $\cE$ be the class of step functions on $\bbR$, i.e. $f\in\cE$ if and only if $f=\sum_{i=1}^K \theta_i 1_{]x_i,x_{i+1}]}$ for some $\theta_1,\cdots,\theta_K$ and $x_1<\cdots< x_{K+1}$. For $f\in \cE$, we define the random variable $W[f]\in L^2(\Omega)$ by
$$W[f]=\sum_{i=1}^K \theta_i (W(x_{i+1})-W(x_i)) $$
and use the notation $W[f]=\int_{\bbsR} f(x)W(dx)$ to emphasize the analogy with integration (even if the path $W$ need not to have bounded variations on compact sets).\\ 
Define the scalar product $<.,.>_W$ on $\cE$ by 
$$<f_1,f_2>_W={\rm Cov}(W[f_1],W[f_2])\quad , f_1,f_2\in\cE.$$
This is indeed a scalar product since the process $W$ is non degenerate. With these notations, $W: \cE\to L^2(\Omega), f\mapsto W[f]$ is a linear isometry. Define the Hilbert space $L^2_W$ as the closure of $\cE$ with respect to the scalar product $<.,.>_W$. Since $L^2(\Omega)$ is complete and $\cE$ is dense in $L^2_W$ (by construction), we can extend the isometry into  $W:L^2_W\to L^2(\Omega)$. \\
\ \\
{\bf Examples: } If $W$ is the Brownian motion, the scalar product is given by
$$<f_1,f_2>_W=\int_{\bbsR} f_1(u)f_2(u)du,\quad f_1,f_2\in\cE  $$
and then $L^2_W=L^2(\bbR)$. Then $W[f]$ is the Wiener integral of $f$ with respect to Gaussian white noise.\\ 
If $W$ is the fractional Brownian motion of index $H\in (1/2,1)$, then
$$<f_1,f_2>_W=H(2H-1)\int_{\bbsR\times\bbsR} |u-v|^{2H-2}f_1(u)f_2(v)dudv,\quad f_1,f_2\in\cE$$
and elements of $L^2_W$ may not be functions but rather distributions of negative order (see \cite{Nua2} or \cite{Pip}). To our purpose, we will not have to consider this distributions space, but it will be enough to remark that $L^1\cap L^2\subset L^2_W$ and that $|\!|f|\!|_{L^2_W}\leq C_2 |\!|f|\!|_{L^1\cap L^2}$ for some constant $C_2$, where $L^1\cap L^2$ is equipped with the norm $|\!|f|\!|_{L^1\cap L^2}=\max(|\!|f|\!|_{L^1},|\!|f|\!|_{ L^2})$.\\
\ \\*
In order to compare integration with respect to $\mu_h$ and to $W$, we have to suppose that $L^2_W$ contains some suitable subspace of locally integrable functions since $\cF_{\mu_h}\subseteq L^1_{loc}$ is the set of functions $f\in L^1_{loc}$ such that the random series (\ref{eq:serie}) defining $\mu_h[f]$ converges in $L^2(\Omega)$. In view of the above examples, we make the following further assumption:
\begin{itemize}
\item{\bf(RS2.d)} $L^1\cap L^2 \subseteq L^2_W $ with continuous injection, i.e. there is some $C_2>0$ such that for any $f\in L^1\cap L^2$,
$$|\!|f|\!|_{L^2_W}\leq C_2 |\!|f|\!|_{L^1\cap L^2}.$$
\end{itemize}

We then prove the convergence of the random measures $\mu_h$ defined by (\ref{eq:defmu}) to the random noise $W$ on the space $L^1\cap L^2\subseteq L^2_W$. 

\begin{pr}\label{pr:RS2}
Suppose that $\xi$ satisfies assumptions ${\bf(RS2)}$.\\
Then for any $h>0$, $L^1\cap L^2\subset \cF_{\mu_h}$ and if $f_n\to f$ in $L^1\cap L^2$ and $h_n\to 0$ as $n\to\infty$,
$$\mu_{h_n}[f_n]\stackrel{\cL}{\Longrightarrow} W[f]=\int_{\bbsR}f(x)W(dx). $$
\end{pr}
\noindent{\bf Remark:} It is worth noting that in the important case when $W$ is a Brownian motion, assumption ${\bf(RS2.d)}$ is trivial. Proposition \ref{pr:RS2} applies under assumptions ${\bf(RS2.b-c)}$ and the conclusion can be strengthened: the convergence holds on the whole space $L^2=L^2_W=\cF_{\mu_h}$ rather than only on the subspace $L^1\cap L^2$.

\noindent{\bf Proof:}\\
For the sake of clarity, we divide the proof into four steps.\\*
\noindent{\bf Step 1:} First note that in the simple case when $f_n\equiv f\in \cE$, Proposition \ref{pr:RS2} directly follows from assumption ${\bf(RS2.b)}$. Let $f=\sum_{i=1}^K \theta_i 1_{]x_i,x_{i+1}]}$. Using equation (\ref{eq:lien}) ,
$$\mu_{h_n}[f]=\sum_{i=1}^K \theta_i (W_{h_n}(x_{i+1})-W_{h_n}(x_i)).$$
Then, assumption $\mbox{\bf(RS2.b)}$ states that as $n\rightarrow +\infty$, 
$$\mu_{h_n}[f]\stackrel{\cL}{\Longrightarrow} \sum_{i=1}^K \theta_i (W(x_{i+1})-W(x_i))=W[f].$$
\noindent{\bf Step 2:} We now prove that for every fixed integer $n$, $L^1\cap L^2\subset \cF_{\mu_{h_n}}$ and that $\mu_{h_n}$ induces a linear functional $L^1\cap L^2\to L^2(\Omega)$ with norm operator bounded by $C_1C_2$.\\ 
The linearity of $\mu_{h_n}$ on $\cE$ is straightforward. Using equation (\ref{eq:lien}), assumption ${\bf(RS2.c)}$ can be rewritten as 
$$\forall f\in\cE,\quad |\!|\mu_{h_n}[f]|\!|_{L^2(\Omega)}\leq C_1 |\!|W[f]|\!|_{L^2(\Omega)}=C_1 |\!|f|\!|_{L^2_W},$$
and then, in view of assumption ${\bf(RS2.d)}$, for every $f\in \cE$,
$$|\!|\mu_{h_n}[f]|\!|_{L^2(\Omega)}\leq C_1 C_2 |\!|f|\!|_{L^1\cap L^2}.$$
Since $\cE$ is dense in $L^1\cap L^2$ and $L^2(\Omega)$ is complete, the linear functional $\mu_{h_n}$ can be extended to $L^1\cap L^2$ and its norm operator is bounded by $C_1C_2$.

\noindent{\bf Step 3:} We consider now the case when $f_n\equiv f$ belongs to $L^1\cap L^2$.\\
We approximate $f$ by a sequence of simple functions $f_m$ such that $f_m\to f$ in $L^1\cap L^2$ as $m\to\infty$.
Let $Y_n=\mu_{h_n}[f]$, $X_{m,n}=\mu_{h_n}[f_m]$, $X_m=W[f_m]$, $X=W[f]$.
We apply Theorem 4.2 in \cite{Bil} by checking that :\\
- for fixed $m$, $X_{m,n}$ weakly converges to $X_m$ (use step 1 with the simple function $f_m$),\\
- $X_m$ weakly converges to $X$ because using assumption ${\bf(RS2.d)}$,  
$$|\!|X_m-X|\!|_{L^2(\Omega)}=|\!|f_m-f|\!|_{L^2_W}\leq |\!|f_m-f|\!|_{L^1\cap L^2}\to 0,$$
- $\lim_{m\infty}\limsup_{n\infty} \bbP\left(|X_{m,n}-Y_n|\geq \varepsilon\right)=0$ for any $\varepsilon >0$. \\
The last condition is satisfied since
\begin{eqnarray*}
\bbP\left(|X_{m,n}-Y_n|\geq \varepsilon\right)&=&\bbP\left(|\mu_{h_n}[f_m-f]|\geq \varepsilon\right)\\
&\leq& \varepsilon^{-2}|\!|\mu_{h_n}[f_m-f]|\!|_{L^2(\Omega)}^2 \leq (C_1C_2)^2\varepsilon^{-2} |\!|f_m-f|\!|_{L^1\cap L^2}^2.
\end{eqnarray*}
Using Theorem 4.2 in \cite{Bil}, we deduce that $Y_n$ weakly converges to $X$ as $n$ tends to infinity.

\noindent{\bf Step 4:} We consider the general case $f_n\to f$ in $L^1\cap L^2$, and we write
$$\mu_{h_n}[f_n]=\mu_{h_n}[f]+\mu_{h_n}[f_n-f].$$
From step 3,  $\mu_{h_n}[f]$ weakly converges to $W[f]$. From step 2, 
$$|\!|\mu_{h_n}[f_n-f]|\!|_{L^2(\Omega)}\leq C_1 C_2 |\!|f_n-f|\!|_{L^1\cap L^2}$$ 
and hence $\mu_{h_n}[f_n-f]$ converges to $0$ in $L^2(\Omega)$. As a consequence, 
$\mu_{h_n}[f_n]$ weakly converges to $W[f].$ \ \ \ $\Box$

\noindent Note that the above proposition can be applied in the following cases.

\noindent{\bf Example: Weakly dependent sceneries.}\\
Under the assumption of a $L^2$ stationary scenery satisfying some $\theta$-mixing condition, a central limit theorem for triangular arrays implying assumption ${\bf(RS2.b)}$ was proved in \cite{GPP1} (Theorem 3.1) as well as a suitable control of the covariance implying assumption ${\bf(RS2.c)}$ (see Lemma 7.1 in \cite{GPP1}). \\*
By applying the techniques developed in \cite{Leb}, a Kesten-Spitzer result for a simple $\bbZ$-random walk in a random scenery coming from a partially hyperbolic dynamical system can be established. It holds when the scenery satisfies some strong decorrelation properties given in Theorem 1 of \cite{PE2}. We just mention that conditions ${\bf(RS2.a-b)}$ should hold under these last assumptions.   

\noindent{\bf Example: Moving averages sceneries.}\\
Stationary sceneries obtained from moving averages of i.i.d. square integrable random variables are considered in \cite{DOM}, i.e.
$$\xi_k=\sum_{i\in\bbsZ}c_{k-i}\eta_i$$
where $(\eta_i)_{i\in\bbsZ}$ is an i.i.d. square integrable sequence. In the case when $c$ is summable, then the convergence of the rescaled cumulative random scenery to a Brownian motion is shown and assumptions ${\bf(RS2.b-c)}$ hold. \\
In the case when $c$ slowly decays, strong correlations persist in the limit. More precisely, if
$$\lim_{k\to +\infty}|k|^{-\gamma}c_k = p_1 \quad \mbox{and}\quad \lim_{k\to -\infty}|k|^{-\gamma}c_k = p_2 $$
with $\gamma=(1/2,1)$ and $p_1p_2\neq 0$, then convergence of the rescaled scenery to a fractional Brownian motion of index $H=3/2-\gamma\in(1/2,1)$ is shown and assumptions ${\bf(RS2.b-c)}$ are satisfied. Note that in \cite{DOM}, the dimension $d$ of the scenery needs not to be equal to one (spatial sceneries) and the innovations $\eta_i$ are more generally assumed to belong to the domain of attraction of a $\beta$-stable distribution with $1<\beta\leq 2$. In this framework, the convergence of the random measures $\mu_h$ to a fractional $\beta$-stable random noise on the space $L^\beta$ is shown. The method of convergence of random measures to a random noise is thus very robust. We will show how this can be used to obtain new results on RWRS.

\section{Limit theorems for RWRS}
The above results can be applied to prove the convergence of RWRS and to obtain limit theorems under fairly general assumptions.
The main idea is that the RWRS can be rewritten as a function of the rescaled local time $L_n(.,.)$ of the random walk $S$ and of the random measures $\mu_h$ associated with the random scenery $\xi$. Using the scaling property (\ref{eq:scaling}) and introducing the suitable scalings for the local time (equation (\ref{eq:LTscaling})) and for the random measures (equation (\ref{eq:defmu})), we get
\begin{eqnarray*}
Z_{[nt]}&=&\sum_{k=0}^{[nt]}\xi_{S_k} \\
&=& \sum_{x\in\bbsZ} N([nt],x)\xi_x \\
&=& \gamma_1^{-1}\int_{\bbsR} N([nt],[x])\mu_1(dx) \\
&=& \frac{n}{\gamma_{a_n^{-1}}} \mu_{a_n^{-1}} [L_n(t,.)]
\end{eqnarray*}
Or equivalently,  
\begin{equation}\label{eq:fondaRWRS}
n^{-1}\gamma_{a_n^{-1}}Z_{[nt]}=\mu_{a_n^{-1}} [L_n(t,.)]
\end{equation}
with $a_n^{-1}\to 0$. We then use the convergence of the local time and of the random measure to prove the convergence of the RWRS.
Two cases are considered and stated in the same following theorem.
\begin{theo}\label{th:RWRS}
Suppose that the random walk $S$ and the random scenery $\xi$ are independent and verify either:
\begin{enumerate}
\item (First case)  $\xi$ satisfies ${\bf(RS1)}$ for some $\beta\in (1,2]$  and  $S$ satisfies ${\bf(RW)}$ with $p=\beta$.
\item (Second case) $\xi$ satisfies ${\bf(RS2)}$ and  $S$ satisfies ${\bf(RW)}$ with $p=2$.
\end{enumerate}
Then, the RWRS satisfy the following convergence in the sense of the finite dimensional distributions:
$$(n^{-1}\gamma_{a_n^{-1}} Z_{[nt]})_{t\geq 0}\stackrel{\cL}{\Longrightarrow} \left(\int_{\bbsR}L(t,x)W(dx)\right)_{t\geq 0} $$
where $(W(x))_{x\in\bbsR}$  and $(L(t,x))_{t\geq 0;x\in\bbsR}$ are defined in ${\bf(RS)}$ and ${\bf(RW)}$ respectively and independent of each other. 
\end{theo}

\noindent{\bf Proof:}\\*
We only prove the theorem in the first case. 
Thanks to Proposition \ref{stalev} and Skorohod's representation theorem, there exists one probability space on which are defined copies of the processes $L_n$ and $L$ denoted by $\tilde{L}_n$ and $\tilde{L}$ such that the sequence 
$(\tilde{L}_n)_n$ converges almost surely in $L^{p}$ to $\tilde{L}$ as $n$ tends to infinity.
Then, by combining equality (\ref{eq:fondaRWRS}) and Proposition \ref{clem}, we get that conditionally on these local times, for every $m\geq 1$, for every $\theta_1\in\bbR,\ldots,\theta_m\in\bbR$, for every $0\leq t_1<\ldots<t_m$,
\begin{equation}\label{eq:fondaRWRS2}
n^{-1}\gamma_{a_n^{-1}}\sum_{i=1}^m \theta_i Z_{[nt_i]}= \sum_{i=1}^m \theta_i \mu_{a_n^{-1}} [\tilde{L}_n(t_i,.)]
\end{equation}
converges in distribution to 
$$\sum_{i=1}^m \theta_i \int_{\bbsR}\tilde{L}(t_i,x)\, W(dx).$$
The result then easily follows using dominated convergence theorem. \ \  $\Box$

We finally state a result for many walkers in random scenery. We consider $c_n$ independent and identically distributed random walks $S^{(i)}, i=1,\ldots,c_n$ defined either as in Section 2.2 or as in Section 2.3. They are assumed to evolve in the same random scenery $\xi$. The $i$-th random walk induces a RWRS denoted by $Z^{(i)}$, but the processes $Z^{(1)},\cdots,Z^{(c_n)}$ are not independent since all random walkers evolve in the same scenery (see \cite{CS,CD,domgui}). The total random reward at time $n$ is then given by 
$$Z_{n,c_n}=\sum_{i=1}^{c_n}Z_n^{(i)}.$$
Using equation (\ref{eq:fondaRWRS}), we get
\begin{equation}
c_n^{-1}n^{-1}\gamma_{a_n^{-1}}Z_{[nt],c_n}=\mu_{a_n^{-1}} [c_n^{-1}\sum_{i=1}^{c_n}L_n^{(i)}(t,.)].
\end{equation}
\noindent
For a finite number $c_n\equiv c$ of independent walkers in a random scenery the functional approach easily gives us the following convergence
$$(n^{-1}\gamma_{a_n^{-1}} \sum_{i=1}^c Z_{[nt]}^{(i)})_{t\geq 0}\stackrel{\cL}{\Longrightarrow} \left(\int_{\bbsR}(\sum_{i=1}^c L^{(i)}(t,x))W(dx)\right)_{t\geq 0}.$$
When $c_n\to\infty$ as $n\to \infty$, thanks to Propositions \ref{manywa} and \ref{clem}, we get the following result, whose proof is very similar to the proof of Theorem \ref{th:RWRS} and is omitted. 
\begin{theo}
Suppose that the random walk $S^{(i)}$ and the random scenery $\xi$ are independent and verify either:
\begin{enumerate}
\item (First case)  $\xi$ satisfies ${\bf(RS1)}$ for some $\beta\in (1,2]$  and  the $S^{(i)}$ verify ${\bf(RW)}$ as well as the integrability condition (\ref{eq:integLT}) with $p=\beta$.
\item (Second case)  $\xi$ satisfies ${\bf(RS2)}$   and  the $S^{(i)}$ verify ${\bf(RW)}$ as well as the integrability condition (\ref{eq:integLT}) with $p=2$.
\end{enumerate}
Then, the total random reward satisfies the following convergence in the sense of the finite dimensional distributions:
$$(c_n^{-1}n^{-1}\gamma_{a_n^{-1}} Z_{[nt],c_n})_{t\geq 0}\stackrel{\cL}{\Longrightarrow} \left(\int_{\bbsR}\bbE(L(t,x))W(dx)\right)_{t\geq 0}.$$
where $(W(x))_{x\in\bbsR}$  and $(L(t,x))_{t\geq 0;x\in\bbsR}$ are defined in ${\bf(RS)}$ and ${\bf(RW)}$ respectively.
\end{theo}

\section{Conclusion and further comments}
In this paper a functional approach for the analysis of RWRS was proposed. The method is quite powerful and allows to recover many of the previous results about RWRS as well as new ones. However, our results are limited to the convergence of the finite-dimensional distributions and the question of the functional convergence (for example in the Skohorod topology) naturally arises. The classical techniques used for proving tightness such as Kolmogorov's criterion can be adapted. The reader should refer to \cite{CD}, \cite{domgui}, \cite{KS} or \cite{WA} where tightness results were proved. In a few words, in the case of a random scenery satisfying ${\bf(RS2)}$, the second order moment of the increments of the RWRS is given for any $0\leq t_1<t_2$ by
\begin{eqnarray*}
\bbE\left[|n^{-1}\gamma_{a_n^{-1}} Z_{[nt_2]}-n^{-1}\gamma_{a_n^{-1}} Z_{[nt_1]}|^2 \right]
&=&\bbE\left[| \mu_{a_n^{-1}} [L_n(t_2,.)-L_n(t_1,.)] |^2 \right]\\
&\leq & (C_1C_2)^2 \bbE\left[\int_{\bbsR} | L_n(t_2,x)-L_n(t_1,x) |^2 dx\right].
\end{eqnarray*} 
The proof of the tightness then relies on the analysis of the self-intersection local time of the random walk (see \cite{KS} for RW with i.i.d. increments or \cite{WA} for strongly correlated Gaussian RW).


\addcontentsline{toc}{section}{\bf References}
\begin{thebibliography}{12}
\bibitem{Bil}{\sc Billingsley, P.} (1968) {\it Convergence of probability measures.} {\rm Wiley, New York.}

\bibitem{bol} {\sc Bolthausen, E.} (1989) {\rm A central limit theorem for two-dimensional random
walks in random sceneries.} {\it Ann. Probab.,} Vol. {\bf 17}, 108 -- 115.

\bibitem{Bor1} {\sc Borodin, A.N.} (1979) {\rm Limit theorems for sums of independent random variables defined on a nonrecurrent random walk.} {\it Investigations in the theory of probability distributions, IV. Zap. Nauchn. Sem. Leningrad. Otdel. Mat. Inst. Steklov.(LOMI)}, Vol. {\bf 85}, 17--29.

\bibitem{Bor2} {\sc Borodin, A.N.} (1983) {\rm Limit theorems for sums of independent random variables defined on a recurrent random walk.} {\it Theory Probab. Appl.}, Vol. {\bf 28}, 98--114.

\bibitem{C} {\sc Cadre, B.} (1995) {\rm Etude de convergence en loi de fonctionnelles de processus: Formes quadratiques ou multilin\'eaires al\'eatoires, Temps locaux d'intersection de marches al\'eatoires, Th\'eor\`eme central limite presque s\^{u}r.} {\it PHD Thesis, Universit\'{e} Rennes 1.}

\bibitem{CHK} {\sc Chen, X.} and {\sc Khoshnevisan, D.} (2009) {\rm From charged polymers to random walk in random scenery.} To appear in {\it Proceedings of the Third Erich L. Lehmann Symposium.}

\bibitem{CH} {\sc Chen, X.} and {\sc Li, W.V.} (2004) {\rm Large and moderate deviations for intersection local times.} {\it Probab. Theory Relat. Fields}, Vol. {\bf 128}, 213--254. 

\bibitem{CD} {\sc Cohen, S.} and {\sc Dombry, C.} (2009) {\rm Convergence of dependent walks in a random scenery to fBm-local time fractional stable motions.} To appear in {\it Journal of Mathematics of Kyoto University}.

\bibitem{CS} {\sc Cohen, S.} and {\sc Samorodnitsky, G.} (2006) {\rm Random rewards, fractional
Brownian local times and stable self-similar processes.} {\it Ann. Appl.
Probab.,} Vol. {\bf 16}, No 3, 1432--1461.

\bibitem{DM} {\sc Dobrushin, R.L.} and {\sc Major, P.} (1979) {\rm Non-central limit theorems for non-linear functionals of Gaussian
fields.} {\it Z. Wahrscheinlichkeitstheorie verw. Gebiete} 50, 27-52.

\bibitem{DOM} {\sc Dombry,C.} Convergence to stable noise and applications. Preprint.

\bibitem{domgui} {\sc Dombry, C.} and {\sc Guillotin-Plantard, N.} (2009) {\rm Discrete approximation
of a stable self-similar stationary increments process.} {\it Bernoulli}, Vol. {\bf 15}, No 1, 195--222. 

\bibitem{GPLN1} {\sc Guillotin-Plantard, N.} and {\sc Le Ny, A.} (2007) {\rm Transient random walks on 2d-oriented lattices.} {\it Theory of Probability and Its Applications (TVP)}, Vol. {\bf 52}, No 4, 815--826.

\bibitem{GPLN2} {\sc Guillotin-Plantard, N.} and {\sc Le Ny, A.} (2008) {\rm A functional limit theorem for a 2d- random walk with dependent marginals.} {\it Electronic Communications in Probability}, Vol. {\bf 13}, 337--351.

\bibitem{GPP1} {\sc Guillotin-Plantard, N.} and {\sc Prieur, C.} (2009) {\rm Central limit theorem for sampled sums of dependent random variables.} To appear in {\it ESAIM P\&S}.

\bibitem{GPP2} {\sc Guillotin-Plantard, N.} and {\sc Prieur, C.} (2009) {\rm Limit theorem for random walk in weakly dependent random scenery.} Submitted.

\bibitem{KS}{\sc Kesten, H.} and {\sc Spitzer, F.}  (1979) {\rm A limit theorem related to a new class of self-similar processes.} {\it Zeitschrift für Wahrcheinlichkeitstheorie und verwandte Gebiete}, Vol. {\bf 50},  5--25.

\bibitem{LX} {\sc Lang, R.} and {\sc Xanh, N.X.} (1983) {\rm Strongly correlated random fields as observed by a random walker.} {\it Probab. Theory and Related Fields}, Vol. {\bf 64}, 327--340. 

\bibitem{Leb} {\sc Le Borgne, S.} (2006) {\rm Exemples de systèmes dynamiques quasi-hyperboliques à décorrélations lentes.} {\it C.R.A.S.}, 343, N. 2, 125--128.

\bibitem{led} {\sc Ledoux, M.} and {\sc Talagrand, M.} (1991) {\rm Probability in Banach spaces. Isoperimetry and processes.} Ergebnisse der Mathematik und ihrer Grenzgebiete (3) [Results in Mathematics and Related Areas (3)], 23. Springer-Verlag, Berlin.

\bibitem{Mae} {\sc Maejima, M.} (1996) {\rm Limit theorems related to a class of operator-self-similar processes.} {\it Nagoya Math. J.} 142, 161--181.

\bibitem{Nua} {\sc Nualart, D.} (2003) {\rm Stochastic integration with respect to fractional Brownian motion and applications.} {\it Contemporary Mathematics} 336, 3--39.

\bibitem{Nua2} {\sc Nualart, D.} (2006) {\rm The Malliavin calculus and related topics.} {\it Probability and its Applications} (New York). Springer-Verlag, Berlin. Second edition.

\bibitem{PE2} {\sc Pène, F.} (2007) {\rm Transient random walk in $Z^2$ with stationary orientations.} To appear in {\it ESAIM P\&S}.

\bibitem{Pip} {\sc Pipiras, V.} and {\sc Taqqu, M.S.}(2000) {\rm Integration questions related to fractional
Brownian motion.} {\it Probab. Theory Rel. Fields} 118,  121--291.

\bibitem{WA} {\sc Wang, W.} (2003) {\rm Weak convergence to fractional Brownian motion in Brownian scenery.}
{\it Probab. Theory and Related Fields}, Vol. {\bf 126}, 203--220.

\end{thebibliography}
\end{document}